\documentclass[11pt]{article}
\usepackage{amsmath,amsthm,amsfonts}
\DeclareMathOperator{\Hom}{Hom}

\newtheorem{thm}{Theorem}[section]
 \newtheorem{lemma}[thm]{Lemma}

 \newtheorem{cor}[thm]{Corollary}
\theoremstyle{definition}
 \newtheorem{defn}[thm]{Definition}
\theoremstyle{definition}

\theoremstyle{definition}
 
\def\End{\mbox{\rm End}\,}

\def\Hom{\mbox{\rm Hom}\,}

\def\Im{\mbox{\rm Im}\,}

\def\id{\mbox{\rm id}}

\def\into{\hookrightarrow}
\def\to{\rightarrow}

\def\o{\otimes}    
\def\bra{\langle}
\def\ket{\rangle}

\def\cop{\Delta}

\def\eps{\varepsilon}

\def\du1{\hat 1}

\def\-1{_{(-1)}}
\def\0{_{(0)}}
\def\1{_{(1)}}
\def\2{_{(2)}}
\def\3{_{(3)}}
\def\4{_{(4)}}

\def\|{\, | \,}

\def\du1{\hat 1}
\def\lact{\triangleright}

\def\op{^{\rm op}}

\begin{document}
\title{\LARGE\bf Simplicial Hochschild
cochains as  \\ an Amitsur complex}
\date{}
\author{\Large Lars Kadison \\ \\
Department of Mathematics, University
of Pennsylvania\\
Philadelphia, PA 19104-6395\\ \\
E-mail: lkadison@math.upenn.edu}
\maketitle
\thispagestyle{empty}
\begin{abstract}
It is demonstrated that the cochain complex
of  relative Hochschild $A$-valued  cochains of
a depth two extension $A \| B$ under cup product
is isomorphic as a differential graded algebra with the Amitsur complex of the coring $S = \End {}_BA_B$  over the centralizer $R = A^B$ with grouplike element $1_S$, which itself is  isomorphic
to the Cartier complex of
$S$ with coefficients in the
$(S,S)$-bicomodule $R^e$.  
This  specializes to finite dimensional algebras, H-separable extensions and Hopf-Galois
extensions.  
\par\smallskip
{\bf 2000 MSC:} 18G25. 
\end{abstract}

\section{Introduction}

Relative Hochschild cohomology of a subring $B \subseteq A$ or ring homomorphism $B \to A$ is set forth in \cite{H}.  The coefficients
of the general form of the  cohomology theory are taken in a bimodule $M$ over $A$.  If $M = A^*$ is the $k$-dual of the $k$-algebra $A$, this gives rise to a cyclic symmetry exploited in cyclic cohomology.  If $M = A$, this has been shown to be related to the simplicial cohomology of a finitely triangulated space via barycentric subdivision, the 
poset algebra of incidence relations and the
separable subalgebra of simplices by Gerstenhaber and Schack in a series of papers beginning with \cite{GS}.  The $A$-valued relative cohomology groups
of $(A,B)$ are also of interest in deformation theory.  
We  refer  to the relative
Hochchild cochains with cohomology groups
$H^n(A,B; A)$ as simplicial Hochschild
cochains with cohomology.  

In this note we will extend the following algebraic  result in \cite{LK2006}: given a depth two ring extension $A \| B$ with centralizer
$R = A^B$ and endomorphism ring $S = \End {}_BA_B$, the simplicial
Hochschild cochains under cup product are
isomorphic as a graded algebra to the tensor algebra  of the $(R,R)$-bimodule
$S$.  Since $S$ is a left bialgebroid over
$R$, it is in particular an $R$-coring with grouplike element $1_S = \id_A$.  The Amitsur complex of such a coring is a differential graded algebra explained in \cite[29.2]{BW}. We note below that the algebra isomorphism in \cite{LK2006}
extends to an isomorphism of differential
graded algebras. We also note that the Amitsur complex of the underlying
coring of a bialgebroid is a Cartier complex with coefficients in a
bicomodule formed from source and target homomorphisms. 
We remark on the consequences for relative Hochschild cohomology of various
types of Galois extensions with bialgebroid
action or coaction.

\section{Preliminaries on depth two extensions}

All rings and algebras are unital associative; homomorphisms and modules are unital as well.
Let $R$ be a ring, and $M_R$, $N_R$ be right $R$-modules.  The notation $M / N$ denotes that $M$ is $R$-module isomorphic to a direct summand of an $n$-fold direct sum power of $N$:  $M \oplus * \cong N^n$.  Recall that $M$ and $N$
are similar \cite[p. 268]{AF} if $M / N$ and $N / M$.
A ring homomorphism $B \to A$ is sometimes called a ring extension $A \| B$ (proper
ring extension if $B \into A$).  

\begin{defn}
A ring homomorphism $B \to A$ is said to be a right
depth two (rD2) extension if the natural $(A,B)$-bimodules $A \o_B A$ and $A$ are
similar.  
\end{defn}
Left D2 extension is defined similarly using the natural $(B,A)$-bimodule structures: a
D2 extension is both rD2 and $\ell$D2.   
Note that in either case any ring extension satisfies $A / A \o_B A$.  

Note some obvious cases of depth two:
1) $A$ a finite dimensional algebra, $B$
the ground field. 
 2) $A \| B$ an H-separable
extension. 3) $A \| B$ a finite Hopf-Galois
extension, since the Galois isomorphism
$A \o_B A \stackrel{\cong}{\longrightarrow}
A \o H$ is an $(A,B)$-bimodule arrow
(and its twist by the antipode shows
$A \| B$ to be $\ell$D2 as well).  

Fix the notation $S := \End {}_BA_B$ and $R = A^B$.
Equip $S$ with $(R,R)$-bimodule structure $$r\cdot \alpha \cdot s
= r\alpha(-)s = \lambda_r \circ \rho_s \circ \alpha $$
where $\lambda, \rho: R \to S$ denote left and right multiplication
of $r,s \in R$ on $A$.

\begin{lemma}\cite[3.11]{KS}
If $A \| B$ is rD2, then the module $S_R$ 
is a projective generator and 
\begin{equation}
f_2: \ S \o_R S \stackrel{\cong}{\longrightarrow} \Hom ({}_BA \o_B A_B, {}_BA_B)
\end{equation}
via $f_2(\alpha \o_R \beta) (x \o_B y) = \alpha(x)\beta(y)$ for $x, y \in A$.
\end{lemma}

For example, if $A$ is a finite dimensional algebra over ground field
$B$, then $S = \End A$, the linear endomorphism algebra.  
If $A \| B$ is H-separable, then $S \cong R \o_Z R\op$,
where $Z$ is the center of $A$ \cite[4.8]{KS}.  If $A \| B$ is an
$H^*$-Hopf-Galois extension, then $S \cong R \# H$, the smash product where $H$ has dual action on $A$ restricted to $R$ \cite[4.9]{KS}.

Recall that a left $R$-bialgebroid $H$ is a type
of bialgebra over a possibly noncommutative
base ring $R$. More specifically, $H$ and $R$ are rings with 
``target'' and ``source'' ring anti-homomorphism
and homomorphism $R \to H$, commuting at all values in $H$,
which induce an $(R,R)$-bimodule structure
on $H$ from the left.  W.r.t.\ this structure, there is an $R$-coring structure
$(H,R, \cop, \eps)$ such that $1_H$ is a
grouplike element (see the next section)
and the left $H$-modules form a tensor
category with fiber functor to the category of $(R,R)$-bimodules.  
One of the main theorems in depth two theory is 
\begin{thm}\cite[3.10, 4.1]{KS}
Suppose $A \| B$ is a left or right D2
ring extension.  Then the endomorphism ring $S := \End {}_BA_B$ is a left bialgebroid over the centralizer $A^B : = R$ via the source map $\lambda: R \into S$,
target map $\rho: R\op \into S$, coproduct
\begin{equation}
\label{eq: cop}
f_2(\cop(\alpha))(x \o_B y) = \sum_{(\alpha)} f_2(\alpha\1 \o_R \alpha\2)(x \o_B y) =  \alpha(xy).
\end{equation}
 Also $A$ under the natural action of $S$ is a left $S$-module algebra
with invariant subring $A^S \cong \End {}_EA$ where $E := \End A_B
\stackrel{\cong }{\longleftarrow} A \# S$ via $a \o_R \alpha \mapsto
\lambda_a \circ \alpha$.  
\end{thm}
We note in passing the measuring axiom of module algebra action from
eq.~(\ref{eq: cop}): in Sweedler notation, $\sum_{(\alpha)} \alpha\1(x)\alpha\2(y) = \alpha(xy)$.  Note too that
$\cop(\lambda_r) = \lambda_r \o 1_S$ and $\cop(\rho_s) = 1_S \o_R \rho_s$
for $r,s \in  R$. 
\section{Amitsur complex of a coring with grouplike}

An $R$-coring $\mathcal{C}$ has coassociative coproduct
$\cop: \mathcal{C} \to \mathcal{C} \o_R \mathcal{C}$ and counit
$\eps: \mathcal{C} \to R$, both mappings being $(R,R)$-bimodule homomorphisms.  We assume that $\mathcal{C}$ also has a grouplike
element $g \in \mathcal{C}$, which means that $\cop(g) = g \o_R g$
and $\eps(g) = 1$.  The Amitsur complex $\Omega(\mathcal{C})$ of $(\mathcal{C}, g)$ has
$n$-cochain modules $$ \Omega^n(\mathcal{C}) = \mathcal{C} \o_R \cdots \o_R \mathcal{C}$$
($n$ times $\mathcal{C}$), the zero'th given by  $\Omega^0(\mathcal{C}) = R$.
The Amitsur complex  is the tensor algebra
$$\Omega(\mathcal{C}) = \oplus_{n=0}^{\infty} \Omega^n(\mathcal{C})$$ with a compatible differential $d = \{d^n \}$ where $d^n: \Omega^n(\mathcal{C}) \to \Omega^{n+1}(\mathcal{C})$.  These are defined by $d^0: R \to \mathcal{C}$,
$d^0(r) = rg - gr$, and 
\begin{eqnarray}
\label{eq: coring differential}
d^n(c^1 \o \cdots \o c^n) &=& g \o c^1 \o \cdots \o c^n + (-1)^{n+1}
c^1 \o \cdots \o c^n \o g \\
 &+& \sum_{i=1}^n (-1)^i c^1 \o \cdots \o c^{i-1}
\o \cop(c^i) \o c^{i+1} \o \cdots \o c^n \nonumber
\end{eqnarray}
Some computations show that $\Omega(\mathcal{C})$ is a differential
graded algebra \cite{BW}, with defining equations, $d \circ d = 0$ as well as
the graded Leibniz equation on homogeneous elements,  $$d(\omega \omega')
= (d\omega)\omega' + (-1)^{|\omega|} \omega d\omega'.$$

The name Amitsur complex comes from the case of a ring homomorphism
$B \to A$ and $A$-coring $\mathcal{C} := A \o_B A$ with coproduct
 $\cop(x \o_B y) = x \o_B 1_A \o_B y$ and counit $\eps(x \o_B y) = xy$.
The element $g = 1 \o_B 1$ is a grouplike element.  
We clearly obtain the classical Amitsur complex, which is acyclic if
$A$ is faithfully flat over $B$. In general, the Amitsur complex
of a Galois $A$-coring $(\mathcal{C}, g)$ is acyclic if $A$ is faithfully
flat over the $g$-coinvariants $B = \{ b \in A \| bg = gb \}$ \cite[29.5]{BW}.

The Amitsur complex of interest to this note is the following derivable
from the left bialgebroid $S = \End {}_BA_B$ of a depth two ring extension $A \| B$ with centralizer $A^B = R$.  The underlying $R$-coring $S$ has grouplike element $1_S = \id_A$, with $(R,R)$-bimodule structure, 
coproduct and counit defined in the previous section.  In Sweedler notation, we may summarize
this as follows:
$$ \Omega(S) = R \, \oplus \, S \, \oplus \, S \o_R S \, \oplus \, 
S \o_R S \o_R S \, \oplus  \cdots $$
$$ d^0(r) = \lambda_r - \rho_r, \  d^1(\alpha) = 1_S \o_R \alpha - \alpha\1 \o_R \alpha\2 + \alpha \o_R 1_S, \ldots $$ 

  It is interesting to remark that this particular Amitsur complex is
naturally isomorphic to a Cartier complex of the $R$-coring $S$ with
coefficients in the $(S,S)$-bicomodule $R^e$ \cite[30.3]{BW}.  The right coaction
is given by $\rho^R(r \o s) = r \o \rho_s$, left coaction by
$\rho^L(r \o s) = \lambda_r \o s$, and we note that $\mbox{\rm Hom}_{R-R}(R^e, \Omega^n(S)) \cong \Omega^n(S)$, the differentials
being preserved by the isomorphism.  

\section{Cup product in simplicial Hochschild cohomology}

Let $A \| B$ be an extension of  $K$-algebras.  We briefly recall the $B$-relative
Hochschild cohomology of $A$ with coefficients in $A$ (for coefficients in a bimodule,
see the source \cite{H}).  The zero'th cochain group $C^0(A,B;A) =
A^B = R$, while the $n$'th cochain group $$C^n(A,B;A) = \mbox{\rm Hom}_{B-B}(A \o_B \cdots \o_B A,A)$$ ($n$ times $A$ in the domain).  In particular, $C^1(A,B;A) =  \End {}_BA_B = S$.  The
coboundary $\delta^n: C^n(A,B;A) \to C^{n+1}(A,B;A)$ is given by
\begin{eqnarray}
\label{eq: Hoch}
(\delta^n f)(a_1 \o \cdots \o a_{n+1}) &=& a_1 f(a_2 \o \cdots \o a_{n+1}) 
+ (-1)^{n+1} f(a_1 \o \cdots \o a_n)a_{n+1}   \nonumber \\
 & + & \sum_{i=1}^n  
(-1)^i f(a_1 \o \cdots \o a_i a_{i+1} \o \cdots \o a_{n+1})  
\end{eqnarray}
and $\delta^0: R \to S$ is given by
$\delta^0(r) = \lambda_r -\rho_r$.  The mappings 
satisfy $\delta^{n+1} \circ \delta^n = 0$ for each $n \geq 0 $. Its cohomology is denoted by
$H^n(A,B;A) = \ker \delta^n / \Im \delta^{n-1}$, and might be referred to as a 
simplicial Hochschild cohomology, since  this cohomology
is isomorphic to simplicial cohomology
if $A$ is the  poset algebra of incidence relation in a finite simplicial complex and $B$ is the separable subalgebra of simplices, where $A$ is 
embeddable in an upper triangular matrix algebra with $B$
the diagonal matrices \cite{GS}. 

The cup product $\cup: C^m(A,B;A) \o_K C^n(A,B;A) \to C^{n+m}(A,B;A)$ makes use of the multiplicative
stucture on $A$ and is given by
\begin{equation}
(f \cup g)(a_1 \o \cdots \o a_{n+m}) = f(a_1 \o \cdots \o a_m)g(a_{m+1} \o \cdots \o a_{n+m})
\end{equation}
which satisfies the equation $\delta^{n+m}(f \cup g) = (\delta^m f) \cup g + (-1)^m f \cup \delta^n g$ \cite{GS}.  Cup product therefore passes to a product on the cohomology.  
  We note that $(C^*(A,B;A), \cup, +, \delta)$ is a differential graded algebra we denote by $C(A,B)$.    

\begin{thm}
Suppose $A \| B$ is a right or left D2 algebra extension.  Then the relative Hochschild $A$-valued cochains $C(A,B)$ is   isomorphic as a differential graded algebra to the Amitsur complex $\Omega(S)$ of the $R$-coring $S$.  
\end{thm}

\begin{proof}
We define a mapping $f$ by $f_0 = \id_R$,
$f_1 = \id_S$, and for $n > 1$,  
\begin{equation}
f_n: \ S \o_R \cdots \o_R S \stackrel{\cong}{\longrightarrow} \mbox{\rm Hom}_{B-B}(A \o_B \cdots \o_B A, A).
\end{equation}
by $f_n(\alpha_1 \o \cdots \o \alpha_n)= \alpha_1 \cup \cdots \cup \alpha_n$.
(Note that $f_2$ is consistent with our notation in section~2.)
We proved by induction on $n$ in \cite[Theorem 5.1]{LK2006} that $f$ is an isomorphism of
graded algebras.  We complete the proof
by noting that $f$ is a cochain morphism, i.e., 
commutes with differentials. 
For $n = 0$, we note that $\delta^0 \circ
f_0 = f_1 \circ d^0$, since $d^0 = \delta^0$.  For $n=1$,
$$ \delta^1(f_1(\alpha))(a_1 \o_B a_2) = a_1\alpha(a_2) - \alpha(a_1 a_2)
+ \alpha(a_1)a_2 $$
$$=  f_2(1_S \o_R \alpha - \alpha\1 \o_R \alpha\2 + \alpha \o_R 1_S)(a_1 \o_B a_2) = f_2(d^1(\alpha))(a_1 \o_B a_2) $$
using eq.~(\ref{eq: cop}). The induction step is carried out in a similar but tedious computation:  this completes the proof that $C(A,B) \cong \Omega(S)$. 
\end{proof} 
\section{Applications of the theorem}
We immediately note that the cohomology rings of the two differential
graded algebras are isomorphic.  
\begin{cor}
Relative $A$-valued Hochschild cohomology is isomorphic to the cohomology
of the $A^B$-coring $S = \End {}_BA_B$:
\begin{equation}
H^*(A,B;A) \cong H^*(\Omega(S), d)
\end{equation}
if $A \| B$ is a left or right depth two extension.  
\end{cor}
For example, a depth two f.g.\ projective extension is separable 
iff its $R$-coring $S$ is coseparable \cite[Theorem 3.1]{LK2005}.
Cartier cohomology of a coseparable coring with any coefficients  vanishes
in positive dimensions \cite[30.4]{BW} as does Hochschild cohomology
of a separable extension \cite{H}. 
But cohomology of the Amitsur complex above
is a particular case of Cartier cohomology
as noted at the end of section~3:
\begin{equation}
H^*(\Omega(S), d) \cong H^*_{\rm Ca}(S,R^e).
\end{equation}

\begin{cor}
If the ring extension $A \| B$ is H-separable and one-sided faithfully flat, then the relative Hochschild cohomology, 
$H^n(A,B; A) $ vanishes in positive
dimensions. 
\end{cor}
\begin{proof}
The extension is necessarily proper by faithful flatness.  
Note that $S \cong R \o_Z R$ is a Galois $R$-coring, since
$\{ r \in R \| r\cdot 1_S = 1_S \cdot r \} = Z$, the center of $A$
and the isomorphism $r \o s \mapsto \lambda_r \circ \rho_s$ is clearly an 
$R$-coring homomorphism.  
Whence
$\Omega(S)$ is acyclic by \cite[29.5]{BW}. 
\end{proof}
This  also follows from proving that an H-separable extensions is  separable.      

The next corollary may be stated more generally for algebras over a base ring which
is hereditary, if the universal coefficient theorem is taken into account.  Let
$K$ be a Hopf algebra. 
\begin{cor}
Suppose $A \| B$ is a finite Hopf-$K^*$-Galois extension
of algebras over a field $k$.
Then relative Hochschild $A$-valued cohomology is isomorphic
to the Cartier  cohomology of the underlying coalgebra $K$ with trivial coefficients: for $n \geq 2$, 
\begin{equation}
H^n(A,B;A ) \cong A^B \o_k H^n_{\rm Ca}(K, k).
\end{equation}
\end{cor}
\begin{proof}
This follows from the determination of $ R \o_k K \cong S$
via $r \o h \mapsto \lambda_r \circ (h \lact \cdot)$, and that
$\cop_S = R \otimes \cop_K$
in \cite{KS}. The relation of action of $K$ on $A$ to coaction $A \to A \o K^*$ is given by $h \lact a = a\0 \bra a\1, h \ket$. The $K$-bicomodule structure on $k$ is given by the unit
$k \to K$.   
Note that $\Omega^n(S) \cong
R \o K \o \cdots \o K$ ($n$ times $K$),
where $d^n = R \o d_c^n$ and $d^n_c$ is
the differential for coalgebra cohomology of $K$ with coefficients in $k$ \cite[30.3]{BW}.   
\end{proof}
For example, a finite dimensional Hopf algebra $K$ is Galois over
$k1_K$ via its coproduct as coaction,
where $K^*$ acts on $K$ via $h^* \rightharpoonup h = h\1 \bra h^*, h\2 \ket$. In this case, relative cohomology recovers absolute cohomology
and the corollary states something well-known in a somewhat different
perspective: for $n \geq 2$, 
$
H^n(K,K) \cong K \o H^n_{\rm Ca}(K^*,k)$
(also, $\cong H^n_{\rm Ca}(K^*,K^*)$). 
\section*{Acknowledgement}
The author thanks the organizers and participants of A.G.M.F.\ in Gothenburg and the Norwegian algebra meeting in Oslo (Nov.\ 1-2, 2007) for the stimulating focus on cohomology.

\end{document}